\input amstex\documentstyle{amsppt}  
\pagewidth{12.5cm}\pageheight{19cm}\magnification\magstep1  
\topmatter
\title On the definition of unipotent representations\endtitle
\author G. Lusztig\endauthor
\address{Department of Mathematics, M.I.T., Cambridge, MA 02139}\endaddress
\thanks{Supported by NSF grant DMS-1855773}\endthanks
\endtopmatter   
\document

\define\Irr{\text{\rm Irr}}

\define\si{\sim}

\define\part{\partial}
\define\emp{\emptyset}

\define\ra{\rangle}

\define\m{\mapsto}
\define\do{\dots}
\define\la{\langle}

\define\ti{\tilde}
\define\nl{\newline}
\redefine\i{^{-1}}

\define\bbq{\bar{\QQ}_l}

\redefine\c{\chi}

\define\e{\epsilon}

\define\r{\rho}

\define\th{\theta}

\redefine\l{\lambda}

\define\kk{\bold k}

\define\CC{\bold C}

\define\FF{\bold F}

\define\QQ{\bold Q}

\define\ZZ{\bold Z}

\define\ct{\Cal T}

\subhead 1\endsubhead
Let $G$ be a reductive connected algebraic group over $\kk$, an algebraic closure of a finite
field $\FF_q$. We assume that $G$ has a fixed $\FF_q$-structure with corresponding Frobenius
map $F:G@>>>G$. Let $\ct_1$ be the set of $F$-stable maximal tori of $G$ and let $\ti\ct_1$
be the set of pairs $(T,\th)$ where $T\in\ct_1$ and $\th\in(T^F)^*$, the group of
characters $T^F@>>>\bbq^*$.
(Here $\bbq$ is an algebraic closure of the field of $l$-adic numbers with $l$ a prime
invertible in $\kk$; we set $T^F=\{t\in T;F(t)=t\}$.) Let $\Irr(G^F)$ be the set of
isomorphism classes of irreducible representations of $G^F=\{g\in G;F(g)=g\}$ over $\bbq$.
In \cite{DL} an equivalence relation $\si$ (geometric conjugacy) on $\ti\ct_1$ 
and a map $E:\Irr(G^F)@>>>\ti\ct_1/\si$ was defined (see no.2).
Here $\ti\ct_1/\si$ is the set of equivalence classes for $\si$.
As shown in \cite{DL}, $E$ can be also viewed as a map from $\Irr(G^F)$ to the set of semisimple
classes in the Langlands dual group of $G$ (either over $\kk$ or over $\CC$).
The definition of $E$
given in \cite{DL} is cohomological. In this paper we will give an alternative
definition of $E$ which is elementary (non-cohomological), assuming that $q$ is large
enough. (See Proposition 4.) However, the fact that this
definition is correct still depends on the cohomological arguments in \cite{DL}.
As an application we give an elementary definition of unipotent representations of $G^F$
(when $q$ is large). 

I have proved Proposition 4 in 1977 (it was part of a talk I gave at the LMS Symposium
on Representations of Lie groups, Oxford, June-July 1977). Recently C. Chan amd M. Oi \cite{Ch}
proved that when $q$ is large, an irreducible
cuspidal representation of $G^F$ of the form $\pm R^\th_T$ can be recovered from
the knowledge of its character at regular semisimple elements. This can be deduced from
proposition 4.

\subhead 2\endsubhead
For any integer $n\ge1$ let $\ct_n$ be the set of $F^n$-stable maximal tori of $G$ and
let $\ti\ct_n$ be the set of pairs $(S,\l)$ where $S\in\ct_n$ and
$\l:S^{F^n}@>>>\bbq^*$
is a character. Any $(T,\th)\in\ct_1$ gives rise to a pair $(T,\th N_{T,n})\in\ti\ct_n$;
here $N_{T,n}:T^{F^n}@>>>T^F$ is the homomorphism given by $t\m tF(t)F^2(t)...F^{n-1}(t)$.
Following \cite{DL, 5.5} we say that $(T,\th)\in\ti\ct_1$ and $(T',\th')\in\ti\ct_1$ are
geometrically conjugate if there exists $n\ge1$ and $g\in G^{F^n}$ such that
$gTg\i=T'$ and $\th'(N_{T',n}(gtg\i))=\th(N_{T,n}(t))$ for all $t\in T^{F^n}$. This is an
equivalence relation on $\ti\ct_1$ denoted by $\si$. Let $\ti\ct_1/\si$ be the set of
equivalence classes.

For any $(T,\th)\in \ti\ct_1$ let $R_T^\th$ be the virtual representation of $G^F$ over $\bbq$
defined in \cite{DL, 1.20, 4.3}. Let $\r\in\Irr(G^F)$. Let $\c_\r:G^F@>>>\bbq$ be the character
of $\r$.
For any $T\in\ct_1$ let $Z_{\r,T}=\{\th\in(T^F)^*;\la\r,R_T^\th\ra\ne0\}$; here
$\la\r,R_T^\th\ra\ne0$ is the multiplicity of $\r$ in $R_T^\th$. 
According to \cite{DL, 7.7},

(a) there exists $T\in\ct_1$ such that $Z_{\r,T}\ne\emp$.
\nl
Moreover, by \cite{DL, 6.3}, 

(b) $\{(T,\th)\in\ti\ct_1;\th\in Z_{\r,T}\}$ is contained in a single geometric conjugacy class;
\nl
we denote it by $E(\r)$. Thus we have a well defined map   $E:\Irr(G^F)@>>>\ti\ct_1/\si$.
We show:

(c) For any $T\in\ct_1$ we have $|Z_{\r,T}|\le|W|$.
\nl
(For any finite set $Y$ we denote by $|Y|$ the cardinal of $Y$.)
We write the elements of $Z_{\r,T}$ in a sequence $\th_1,\th_2,\do,\th_r$.
We can find $n\ge1$ such that the sequence $\th_1N_{T,n},\th_2N_{T,n},\do,\th_rN_{T,n}$
in $(T^{F^n})^*$ is contained in a single orbit
of the obvious action of $(NT\cap G^{F^n})/T^{F^n}$ on $(T^{F^n})^*$. (Here $NT$ is the
normalizer of $T$ in $G$.) Since $N_{T,n}$ is surjective, we have
$\th_iN_{T,n}\ne\th_jN_{T,n}$ for $i\ne j$. It follows that
$r\le|(NT\cap G^{F^n})/T^{F^n}|\le|NT/T|=|W|$. This proves (b).

\subhead 3\endsubhead
For any $T\in\ct_1$ let $T^F_{rs}$ be the set of elements of $T^F$ which are regular
semisimple in $G$.
By \cite{DL, 7.6}, for $T\in\ct_1,s\in T^F_{rs}$ and $\r\in\Irr(G^F)$ we have
$$\c_\r(s)=\sum_{\th\in Z_{\r,T}}\th(s)\la\r,R_T^\th\ra.\tag a$$
Let $W$ be the Weyl group of $G$.
If $q$ is large enough compared to the rank of $G$ we clearly have

(b) $|T^F-T^F_{rs}|/|T^F|<2^{-2|W|+1}$ for any $T\in\ct_1$.

\proclaim{Proposition 4} Assume that $q$ satisfies 3(b). For any $\r\in\Irr(G^F)$
there is a unique element $\e\in\ti\ct_1/\si$ such that the following holds:
for any $T\in\ct_1$ there exist $m\in\{0,1,2,\do,|W|\}$, $\th_1,\th_2,\do,\th_m$
distinct in $(T^F)^*$ and $c_1,c_2,\do,c_m$ in $\ZZ-\{0\}$ such that
$(T,\th_i)\in \e$ for $i=1,\do,m$ and such that for any $s\in T^F_{rs}$ we have
$\c_\r(s)=c_1\th_1(s)+c_2\th_2(s)+\do+c_m\th_m(s)$.
\endproclaim
The element $\e=E(\r)$ in no.2 satisfies the requirements by
3(a),2(b),2(c). This proves the existence part. We now prove the uniqueness.
Assume that besides $\e$ we have another element $\e'\in\ti\ct_1/\si$, $\e'\ne\e$
such that $\e'$ satisfies a condition similar to that satisfied by $\e$.
By 3(a),2(a),2(b),2(c),
there exist $T\in\ct_1$, $m\in\{1,2,\do,|W|\}$, $\th_1,\th_2,\do,\th_m$
distinct in $(T^F)^*$ and $c_1,c_2,\do,c_m$ in $\ZZ-\{0\}$ such that
$(T,\th_i)\in \e$ for $i=1,\do,m$ and such that for any $s\in T^F_{rs}$ we have
$\c_\r(s)=c_1\th_1(s)+c_2\th_2(s)+\do+c_m\th_m(s)$.
By our assumption on $\e'$
there exist $m'$ in $\{0,1,2,\do,|W|\}$, $\th'_1,\th'_2,\do,\th'_{m'}$ distinct in $(T^F)^*$
and $c'_1,c'_2,\do,c'_{m'}$ in $\ZZ-\{0\}$ such that
$(T,\th'_i)\in\e'$ for $i=1,\do,m'$ and such that
for any $s\in T^F_{rs}$ we have
$\c_\r(s)=c'_1\th'_1(s)+c'_2\th'_2(s)+\do+c'_{m'}\th'_{m'}(s).$
Since $\e\ne\e'$ we have $\{\th_1,\th_2,\do,\th_m\}\cap\{\th'_1,\th'_2,\do,\th'_{m'}\}=\emp$
so that $\{\th_1,\th_2,\do,\th_m,\th'_1,\th'_2,\do,\th'_{m'}\}$ consists of $m+m'\ge1$
distinct characters satisfying
$c_1\th_1(s)+c_2\th_2(s)+\do+c_m\th_m(s)
=c'_1\th'_1(s)+c'_2\th'_2(s)+\do+c'_{m'}\th'_{m'}(s)$
for all $s\in T^F_{rs}$. Thus, the $m+m'$ characters
$\{\th_1,\th_2,\do,\th_m,\th'_1,\th'_2,\do,\th'_{m'}\}$ are linearly dependent when viewed
as functions $T^F_{rs}@>>>\bbq$. But using 3(b) and \cite{L90, Lemma 8.1} (a refinement of
Dedekind's theorem) we see that
$\{\th_1,\th_2,\do,\th_m,\th'_1,\th'_2,\do,\th'_{m'}\}$ are linearly independent when viewed
as functions $T^F_{rs}@>>>\bbq$. (We use that $m+m'\le2|W|$.) This contradiction shows
that uniqueness holds. The proposition is proved.

\subhead 5\endsubhead
We have thus desribed the map $\r\m E(\r)$ in an elementary way as the map $\r\m\e$
in the proposition (for large $q$).

\subhead 6\endsubhead
Following \cite{DL},
we say that $\r\in\Irr(G^F)$ is unipotent if $1\in Z_{\r,T}$ for
some $T\in\ct_1$. From no.4 we obtain the following result (with $q$ large
enough).

(a) $\r$ is unipotent if and only if for any $T\in\ct_1$, $\c_\r|_{T^F_{rs}}$ is constant.

\subhead 7\endsubhead
In \cite{L84, 4.23} to each $\r\in\Irr(G^F)$ we attach a certain family in a Weyl group
defined in terms of $E(\r)$ (assuming that the centre of $G$ is connected).
The method of proof of Proposition 4 can be used to show that this family can also be
recovered in terms of the character of $\r$ at regular semisimple elements (if $q$ is
large) so that it can be defined in a non-cohomological way.

\enddocument